\documentclass{article}
\usepackage{mathrsfs}
\usepackage{amscd}
\usepackage{amssymb, amsmath, amsthm, graphics}
\usepackage{amsmath,amsfonts,amssymb}
\usepackage{latexsym}
\usepackage{graphicx,psfrag}
\usepackage[all]{xy}
\usepackage{hyperref}

\newtheorem{theorem}{Theorem}[section]
\newtheorem{lemma}[theorem]{Lemma}
\newtheorem{prop}[theorem]{Proposition}
\newtheorem{cor}[theorem]{Corollary}

\theoremstyle{definition}
\newtheorem{definition}[theorem]{Definition}
\newtheorem{example}[theorem]{Example}

\theoremstyle{remark}
\newtheorem{remark}[theorem]{Remark}

\numberwithin{equation}{section}

\newcommand{\RR}{{\mathbb R}}

\newcommand{\CC}{{\mathbb C}}
\newcommand{\ZZ}{{\mathbb Z}}

\newcommand{\SL}{{\rm SL}}
\newcommand{\SO}{{\rm SO}}

\newcommand{\id}{{\rm id}}
\newcommand{\Aut}{{\rm Aut}}
\newcommand{\Diff}{{\rm Diff}}
\newcommand{\ud}{{\rm d}}

\begin{document}

\title{Extending $T^p$ automorphisms over $\RR^{p+2}$ and realizing DE attractors}
\author{Fan Ding, Yi Liu, Shicheng Wang, Jiangang Yao}
\maketitle

\begin{abstract}
In this paper we consider the realization of DE attractors by self-diffeomorphisms
of manifolds. For any expanding self-map $\phi:M\to M$ of a connected,
closed $p$-dimensional manifold
$M$, one can always realize a $(p,q)$-type attractor derived from $\phi$ by a compactly-supported
self-diffeomorphsm of $\RR^{p+q}$, as long as $q\geq p+1$. Thus
lower codimensional realizations are more interesting, related to the knotting problem
below the stable range. We show that for any expanding self-map $\phi$ of a standard
smooth $p$-dimensional
torus $T^p$, there is compactly-supported
 self-diffeomorphism of $\RR^{p+2}$ realizing an attractor derived from $\phi$. A
key ingredient of the construction is to understand automorphisms of
$T^p$ which extend over $\RR^{p+2}$ as a self-diffeomorphism via the
standard unknotted embedding $\imath_p:T^p\hookrightarrow\RR^{p+2}$.
We show that these automorphisms form a subgroup $E_{\imath_p}$ of
$\Aut(T^p)$ of index at most $2^p-1$.
\end{abstract}

2000 AMS subject class 37D45 57R25 57R40

\section{Introduction}

Hyperbolic attractors derived from expanding maps were introduced
by Steve Smale in his celebrated paper \cite{Sm} in the 1960s. Smale
posed four families of basic sets for his Spectral Decomposition
Theorem for the non-wandering set of
self-diffeomorphisms of smooth manifolds: \emph{Group 0} which are zero
dimensional ones
 such as isolated points and the Smale horseshoe; \emph{Group A} and \emph{Group DA},
 both of which are derived from Anosov maps; and
  \emph{Group DE} which are attractors derived from expanding maps.
While the first three families arise easily or automatically from
self-diffeomorphisms of manifolds, it is less obvious whether and
how attractors of Group DE could be realized via self-diffeomorphisms of
manifolds. In this paper, we study the realization problem
of DE attractors. We shall stay in the smooth category.

\begin{definition}\label{ehr}
(1) Let $M$ be an connected, closed $p$-dimensional smooth manifold. A smooth map
$\phi:M\to M$ is said to be \emph{expanding} if for some complete
Riemannian metric on $M$, there exist constants $c>0$, $\lambda>1$
such that for any $x\in N$ and $v\in T_xM$, $\|\ud\phi^m(v)\|\geq c\lambda^m\,\|v\|$,
for every integer $m>0$.

(2) Let $M$, $\phi$ be as above. We say $\phi$ \emph{lifts} to
a hyperbolic bundle embedding $e:E\hookrightarrow E$, if $\pi:E\to M$ is a
compact unit (Euclidean) $q$-dimensional disk bundle over $M$, and $e$ is a smooth self-embedding
of $E$ which descends to $\phi$
under $\pi$, and for some positive number $r<1$,
$e$ sends every fiber $q$-disk to an embedded
$q$-subdisk of radius $r$ in the interior of the target fiber. Suppose there exists such a lift, then
 $\Lambda = \bigcap_{\ i\geq0}\, e^i(E)$ is called a \emph{$(p,q)$-type attractor} derived
from the expanding map $\phi$ (w.r.t. $e$), or simply a \emph{DE attractor}.

(3) Let $M,\phi,E,e$ be as above.
Let $X$ be a $(p+q)$-dimensional smooth manifold, and
 $f: X\to X$ be a self-diffeomorphism of $X$. If there
is a smooth embedding $E\subset X$ such that
$\Lambda=\bigcap_{i\geq0}\, f^i(E)$, we say that $f$ \emph{realizes}
the DE attractor $\Lambda$ in manifolds.
\end{definition}

Expanding maps of closed manifolds are topologically conjugate to expanding infranil
endomorphisms of infranil manifolds (\cite{Gr}). It is also known that any flat manifold
admits an expanding map (\cite{ES}).
In his original paper, Smale only considered trivial $(p+1)$-disk bundle embeddings
which lifts an expanding map $\phi$, but it seems to be reasonable and necessary to allow twisted bundles
here to ensure Theorem \ref{main-realize-stable}. For example, if $M$ is an orientable closed non-spin
flat manifold (e.g. \cite{DSS}) while $X$ is spin, it is impossible
to realize any DE attractor of $M$ on $X$ if we restrict ourselves to trivial bundle lifts, simply because the
normal bundle of any embedding of $M$ into $X$ is also non-spin hence nontrivial.

\begin{prop}\label{main-realize-stable} Suppose $M$ is a connected, closed
$p$-dimensional smooth manifold, $p\geq 1$, and $\phi:M\to M$ is an
expanding map. For any $q\geq p+1$, there is a compactly-supported
self-diffeomorphism of $\RR^{p+q}$, which realizes a $(p,q)$-type
attractor derived from $\phi$.\end{prop}

\begin{remark} A self-diffeomorphism of a smooth manifold is said to be compactly-supported if it
fixes every point outside a compact set. This clear implies that the result holds for any $(p+q)$-dimensional
manifold $X$ besides $\RR^{p+q}$.
\end{remark}

Proposition \ref{main-realize-stable} suggests that the realization
problem of DE attractors of codimensions below the `stable range' is
more interesting. This is basically because in lower codimensions,
for an arbitrary embedding $E\subset X$ (with notations in
Definition \ref{ehr}), if any, the cores of $e^i(E)\subset X$,
for different $i\geq 0$ may be knotted in different ways. In this sense, the
realization problem of DE attractors is essentially about the
knotting problem of embeddings and satellite knot constructions.
On the other hand, as any $(p,q)$-type attractor $\Lambda$ derived from
an expanding map $\phi:M\to M$ is homeomorphic to the inverse limit
$\varprojlim\,(M,\phi)=\left\{(x_0,x_1,\cdots)\in\prod_{n=0}^{\infty}M\,|\,
x_i=\phi(x_{i+1}),\,0\leq i<\infty\right\}$ (sometimes known as a $p$-dimensional \emph{solenoid}),
it cannot be embedded into any $(p+1)$-dimensional closed orientable manifold 
if $M$ is also orientable (\cite{JWZ}). Thus
it becomes natural to wonder whether there are DE attractors
realizable in orientable closed manifolds of codimension $2$.
When $M$ is diffeomorphic to the flat $p$-dimensional torus
$T^p=S^1\times\cdots\times S^1$ ($p$ copies), probably the simplest 
$p$-dimensional manifold admitting expanding maps, the study of
unknotted embeddings of $T^p$ into $\RR^{p+2}$ allows us to give a positive answer:

\begin{theorem}\label{main-realize-T^p} For any expanding
map $\phi: T^p\to T^p$, $p\geq1$,  there is a compactly-supported self-diffeomorphism
of $\RR^{p+2}$ which realizes a $(p,2)$-type attractor derived from
$\phi$.\end{theorem}

To get some idea of the proof of Theorem \ref{main-realize-T^p}, consider the toy case when $p=1$.
An expanding map $\phi:S^1\to S^1$ is nothing but a $d$-fold covering where $d>1$. If $e:E\hookrightarrow
E$ is a disk bundle embedding lifting $\phi$, and $\jmath:E\hookrightarrow \RR^3$ is an embedding, the core of $\jmath(E)$ is a knot $K
\subset\RR^3$, and the core of $\jmath\circ e(E)$ is a satelite knot $K'$ with the companion $K$ and the pattern a braid
in the solid torus of winding number $d$. Note also that $E$ must be a trivial
bundle $S^1\times D^2$ in this case. If there exists a self-diffeomorphism $f$ of $\RR^3$ so that $f\circ\jmath=\jmath\circ e$, $K$, $K'$ should
have homeomorphic exteriors, so
$K$ has to be trivial and the braid could be, for example, a $(d,1)$-cable (in the usual longitude-meridian notations of classical
knot theory). This suggests in the general case we should also consider unknotted, framing untwisted embeddings $\jmath:T^p\times
D^2\hookrightarrow\RR^{p+2}$ (Definitions \ref{unknotted} and \ref{framing}) and bundle embeddings $e:T^p\times D^2\hookrightarrow T^p\times D^2$ respecting the framing.

The real issue when $p\geq 2$ is that $\jmath\circ e$ and $\jmath$ may still be 
non-isotopic even if they have isotopic images, no matter how we choose $e$ and $\jmath$. This is essentially because there
are self-diffeomorphisms of $T^p$ that cannot be extended as a self-diffeomorphism of $\RR^{p+2}$
via a given embedding $\imath:T^p\hookrightarrow\RR^{p+2}$, due to certain
spin obstructions, (\cite{DLWY}). To overcome this difficulty, we say
 two embeddings
$\jmath,\jmath':T^p\times D^2\hookrightarrow\RR^{p+2}$ have different types if they are isotopic up to image but not
isotopic as embeddings, (Definitions \ref{framing}, \ref{types}, cf.
Remark \ref{rmkknot}). We
will show that there are at most finitely many types arising in
$\jmath\circ e^i$ for $i\geq 0$, so that some power of $e$ extends
as a (compactly-supported) self-diffeomorphism of $\RR^{p+2}$ over some
$\jmath\circ e^k$.  The key ingredient is the following theorem,
where $\Aut(T^p)\cong\SL(p,\ZZ)$ is the group of automorphisms on
$T^p$, and $E_{\imath}\leq \Aut(T^p)$ consists of elements which
extend as compactly-supported self-diffeomorphisms of $\RR^{p+2}$
via $\imath$, (cf. Section \ref{Sec-extension} and Definition
\ref{esg}):

\begin{theorem}\label{main-extend}
For the standard unknotted embedding $\imath_p:T^p\hookrightarrow\RR^{p+2}$,
($p\geq1$), the subgroup $E_{\imath_p}$
contains the stabilizer of some nontrivial element in $H^1(T^p;\ZZ_2)$
under the natural action of $\Aut(T^p)$. Hence the index of $E_{\imath_p}$ in $\Aut(T^p)$
is at most $2^p-1$.
\end{theorem}
\begin{remark}
In fact, the index of $E_{\imath_p}$ in $\Aut(T^p)$ is exactly $2^p-1$, combined with the
inequality in the other direction as shown in \cite{DLWY}. Theorem \ref{main-extend} may be
rephrased as there are at most $2^p-1$ modular types (Definition \ref{types}) of unknotted embeddings
(Corollary \ref{modTypeNumber}).
\end{remark}

Our strategy to prove Theorem \ref{main-realize-T^p} is to construct a `favorite' lifting
$e$ of a given expanding map $\phi$ of $T^p$ (represented by a integral matrix 
of expanding eigenvalues in this case) using the Smith normal form of integral
matrices, so that whenever an embedding $\jmath:T^p\times D^2\hookrightarrow
\RR^{p+2}$ is framing untwisted, the same holds for $\jmath\circ e$. Choosing $\jmath$
to be standardly unknotted and framing untwisted, 
we show $\jmath\circ e^i$ are all unknotted of modular types for $i\geq 0$.
Unlike the toy case, this is much less obvious in higher dimensions, and the proof
involves some technical manipulations of integral matrices. After all
these are done, one can apply the finiteness result about modular types to derive Theorem \ref{main-realize-T^p}.   

In Section
\ref{Sec-realizeStable}, we prove Proposition
\ref{main-realize-stable} using a classical unknotting theorem. In
Section \ref{Sec-extension}, we define unknotted embeddings
$\imath:T^p\hookrightarrow\RR^{p+2}$ and their types, and prove
Theorem \ref{main-extend}. In Section \ref{Sec-realization}, we
construct disk bundle embeddings lifting expanding maps, and realize
the corresponding DE attractor by compactly-supported
self-diffeomorphisms of $\RR^{p+2}$, proving Theorem
\ref{main-realize-T^p}. We also define unknotted, framing untwisted
embeddings $\jmath:T^p\times D^2\hookrightarrow\RR^{p+2}$ in Section
\ref{Sec-realization}.

\bigskip\noindent\textbf{Acknowledgement}. The first and the third authors are
partially supported by grant No.10631060 of the National Natural
Science Foundation of China.\\

\section{Realizing DE attractors in large codimensions}\label{Sec-realizeStable}

In this section, we prove Proposition \ref{main-realize-stable}.
This is a consequence of an unknotting theorem of Wen-Ts\"{u}n Wu in the 1950's  
about smooth embeddings into Euclidean spaces of codimension right below the stable range,
(\cite{Wu}, cf. also \cite{Ha} for generalizations).

\begin{proof}[{Proof of Proposition \ref{main-realize-stable}}]
When $p=1$, $M$ is diffeomorphic to $S^1$ and any expanding map $\phi:S^1\to S^1$ is a topologically conjugate
to a non-zero degree covering. The realization of $(1,2)$-attractors, i.e. classical solenoids,
is somewhat well-known, for example, as implicitly contained in \cite{Bo} and
\cite{JNW}. We may assume $p\geq 2$ from now on.

Since $q\geq p+1$, we may pick a Whitney embedding of $\imath:M\hookrightarrow\RR^{p+q}$
of $M$ into $\RR^{p+q}$. By \cite{Gr}, we may assume the expanding map
$\phi:M\to M$ is a covering map induced by an infranil endomorphism.
Now $\phi$ induces an immersion $\imath\circ\phi:M\looparrowright\RR^{p+q}$,
which can be pertubed along normal directions into a smooth embedding $\hat\imath:M\hookrightarrow\RR^{p+q}$ by a
Whitney type argument, such that the image $\hat\imath(M)$
remains in the interior of a compact tubular neighborhood $\mathcal{N}\subset\RR^{p+q}$ of $\imath(M)$.

The result in \cite{Wu} says that for any
connected, closed $p$-dimensional smooth manifold $M$ ($p\geq 2$),
any two smooth embeddings of $M$ into $\RR^{p+q}$ are smoothly
isotopic to each other. Note the connectedness is an indispensable 
assumption here. Applying this unknotting theorem, there is a smooth isotopy $F:M\times
[0,1]\to\RR^{p+q}$, such that $F|_{M\times\{0\}}=\imath$, and
$F|_{M\times\{1\}}=\hat\imath$. By the isotopy extension theorem
(cf. \cite[Theorem 1.3]{Hi}), $F$ extends as a diffeotopy 
$F:\RR^{p+q}\times[0,1]\to\RR^{p+q}$ with compact support. Denote
$f=F|_{\RR^{p+q}\times\{1\}}$. By shrinking the radius of
$f(\mathcal{N})$ by a diffeotopy of $\RR^{p+q}$ supported near $f(\mathcal{N})$ if necessary, we may
assume $f(\mathcal{N})$ is contained in the interior of
$\mathcal{N}$. Identify $\mathcal{N}$ as the unit disk bundle of the
normal bundle $N_{\RR^{p+q}}(\imath(M))$ of $\imath(M)$ (w.r.t. any
Euclidean fiber metric) via a diffeomorphism. By a standard
differential topology argument, we may further assume $f$ is
diffeotoped supported near $f(\mathcal{N})$ so that it maps every fiber disk into a fiber disk of
$\mathcal{N}$, and that it satisfies Definition \ref{ehr} (2). Note
$f:\RR^{p+q}\to\RR^{p+q}$ is also compactly supported by the
construction.

Now $\mathcal{N}$ may be regarded as a $q$-disk bundle over $M$,
and $f|_{\mathcal{N}}:\mathcal{N}\to\mathcal{N}$ may be regarded as a hyperbolic bundle embedding lifted from $\phi$,
via the natural identification. Then clearly $f$ realizes a $(p,q)$-type attractor derived from $\phi$ via the inclusion
$\mathcal{N}\subset\RR^{p+q}$.
\end{proof}

It is remarkable at this point that while Proposition
\ref{main-realize-stable} allows one to realize DE attractors in
fairly arbitrary manifolds of sufficiently large dimensions, it is
still a `local' realization anyways. In contrast, for example, if a
self-diffeomorphism $f$ of a closed, orientable $n$-dimensional
smooth manifold $X$ satisfies that the non-wandering set $\Omega(f)$
is a union of finitely many DE attractors and repellers (i.e.
attractors of $f^{-1}$), then $X$ must be a rational homology
sphere, and each DE attractor/repeller must have codimension $2$,
namely of type $(n-2, 2)$, (\cite{DPWY}, cf. \cite{JNW} for an
example of $n=3$).

\section{Unknotted $T^p$ in $\RR^{p+2}$ and extendable
automorphisms}\label{Sec-extension}

In this section, we introduce and study unknotted embeddings of
$T^p$ into $\RR^{p+2}$, and prove Theorem \ref{main-extend}.

Regard $T^p$ as the standard $p$-dimensional torus obtained by quotienting
$\RR^p$ by its integral lattice. The natural action of $\SL(p,\ZZ)$ on
$\RR^p$ descends to an action on $T^p$, so there is a subgroup $\Aut(T^p)$ of the
orientation-preserving diffeomorphism group $\Diff_+(T^p)$ consisting of the transformations
induced by the action. The elements of $\Aut(T^p)$ will be refered as
\emph{automorphisms} of $T^p$. After choosing a product structure
of $T^p\cong S^1_1\times\cdots\times S^1_p$, one may naturally identify $\Aut(T^p)$
with $\SL(p,\ZZ)$.

We start by investigating some important aspects of unknotted
embeddings. It is reasonable to expect that such embeddings are
fairly simple and symmetric, largely agreeing with our low-dimension
intuition. We will
parametrize $S^1$ and $D^2$ as the unit circle and the compact unit disk of $\CC$, respectively. The
real and imaginary part of $z\in\CC$ are often written as $z_x,z_y$.
The standard basis of $\RR^n$ is
$(\vec{\varepsilon}_1,\cdots,\vec{\varepsilon}_n)$, and the
$m$-subspace spanned by
$(\vec{\varepsilon}_{i_1},\cdots,\vec{\varepsilon}_{i_m})$ will be
written as $\RR^m_{i_1,\cdots,i_m}$. Note there is a natural
inclusion of $\RR^n=\RR^n_{1,\cdots,n}$ into $\RR^{n+1}$.

\begin{figure}[htb]
\centering
\psfrag{a}[]{\scriptsize{$\mathbb{R}^{p-1}_{3,\cdots,{p+1}}$}}
\psfrag{b}[]{\scriptsize{$\vec{\varepsilon}_2$}}
\psfrag{c}[]{\scriptsize{$\vec{\varepsilon}_{p+2}$}}
\psfrag{d}[]{\scriptsize{$\imath_{p-1}(T^{p-1})$}}
\psfrag{e}[]{\scriptsize{$\frac12\cdot\vec{\varepsilon}_2+\frac{1}{4}\cdot$\scriptsize{$\imath_{p-1}(T^{p-1})$}}}
\psfrag{f}[]{\scriptsize{$\imath_p(T^p)$}}
\psfrag{g}[]{\scriptsize{$S^1_p$}}
\includegraphics[scale=.7]{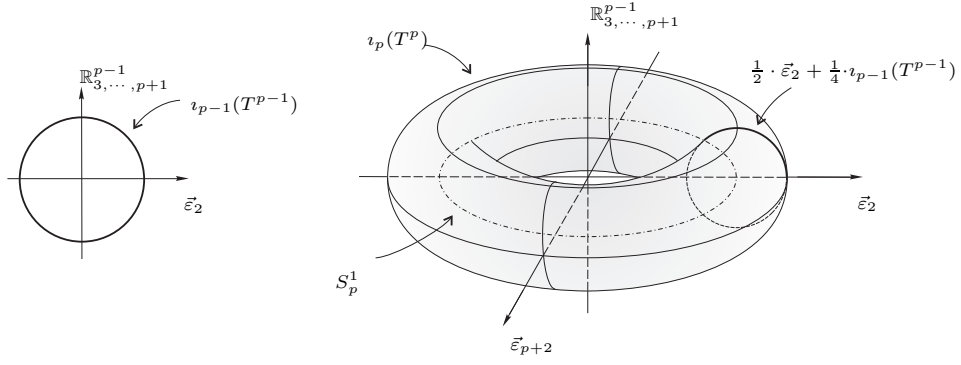}
\caption{The standard model.}\label{figStdModel}
\end{figure}

\begin{example}[The standard model]\label{stdmodel}
Let $\imath_0:{\rm pt}=T^0\hookrightarrow\RR^2$ be the map from the single point to the origin
of $\RR^2$ by
convention. Inductively suppose $\imath_{p-1}$ has been constructed
($p\geq1$) such that $\imath_{p-1}(T^{p-1})\subset {\rm
Int}(D^p)\subset\RR^p_{2,\cdots,p+1}$. Denote
the rotation of $\RR^{p+2}$ on the subspace $\RR^2_{2,p+2}$
 of angle $\arg(u)$ as $\rho_p(u)\in\SO(p+2)$, for any $u\in S^1$. We define
 $\imath_p:T^p=T^{p-1}\times S^1_p$
 as: $$\imath_p(v,u)=\rho_p(u)(\frac12\cdot\vec{\varepsilon}_2+\frac{1}{4}\cdot\imath_{p-1}(v)).$$

This explicitly describes an embedding of
$T^p=S^1_1\times\cdots\times S^1_p$ into
$\RR^{p+1}_{2,\cdots,p+2}$. In Figure \ref{figStdModel}, the images of
$\imath_{p-1}$ and $\imath_p$
are schematically presented on the left and the right respectively. One may imagine
$\vec{\varepsilon}_1$ points perpendicularly outward the page. Observe that the image of $T^p$ is
invariant under $\rho_p(u)$.
\end{example}

The basic feature of the standard model is that it is highly `compressible', in
the sense of the following lemma. Let $B^4\subset\RR^4$ be the compact
$4$-dimensional disk
centered at the origin with radius $2$.

\begin{lemma}\label{k1i} In the standard model for $p\geq2$,
 for each
$i=2,\cdots,p$, the embedding
$\imath_p:T^p=S^1_1\times\cdots\times S^1_p\hookrightarrow \RR^{p+2}$ extends
as an embedding:
$$k_{1i}:B^4\times (S^1_2\times\cdots\times\hat{S}^1_i\times\cdots\times
S^1_{p})\hookrightarrow \RR^{p+2},$$ where $S^1_1\times\cdots\times
S^1_p\subset B^4\times(S^1_2\times\cdots\times
\hat{S^1_i}\times\cdots\times S^1_{p})$ is the standard embedding
(i.e. inclusion)  for $S^1_1\times S^1_i\subset B^4$, and the
identity on other factors.
\end{lemma}

\begin{proof} To see the idea consider a standard $T^2=S^1_1\times S^1_2$
 in $\RR^4$. Make a solid torus $D^2_1\times S^1_2$ by filling
up the $S^1_1$ factor, and attach a semi-sphere in $\RR^4$ along the
core of that solid torus. The result is a `hat' whose regular neighborhood is
diffeomorphic to $D^4$. When $i>2$, one only cares about $S^1_1$ and
$S^1_i$.

The construction is as follows. We first
extend $\imath_1$ as $j_1:D^2_1\hookrightarrow\CC\cong\RR^2_{2,3}$ in an obious fashion. Inductively,
suppose 
$j_{s-1}:D^2_1\times S^1_2\times\cdots\times S^1_{s-1}\hookrightarrow
{\rm Int}(D^s)\subset\RR^s_{2,\cdots,s+1}$ has been constructed. Then
define $j_{s}:(D^2_1\times S^1_2\times\cdots\times
S^1_{s-1})\times S^1_s\hookrightarrow {\rm
Int}(D^{s+1})\subset\RR^{s+1}_{2,\cdots,s+2}$ as:
$$j_s(v,u)=\rho_s(u)(\frac12\cdot\vec{\varepsilon}_2+\frac{1}{4}\cdot j_{s-1}(v)).$$

After $i-1$ steps we obtain $j_i:D^2_1\times
S^1_2\times\cdots\times S^1_i\hookrightarrow {\rm
Int}(D^{i+1})\subset\RR^{i+1}_{2,\cdots,i+2}$. Now let
$\zeta_s(re^{i\theta})$ ($0\leq r\leq 1$) be the rotation of
$\RR^{s+2}$ of angle $\arccos(r)$, on the subspace spanned by
$\rho_s(e^{i\theta})(\vec{\varepsilon}_2)$ and
$\vec{\varepsilon}_1$ (from the former toward the latter). We may
further define $k_{1i,i}:(D^2_1\times S^1_2\times\cdots\times
S^1_{i-1})\times D^2_i\hookrightarrow {\rm Int}(D^{i+2})\subset\RR^{i+2}$ as,
for example,
$$k_{1i,i}(v,re^{i\theta})=\zeta_i(\frac{2r}{1+r^2}e^{i\theta})(j_i(v,e^{i\theta})).$$

Then repeat the standard construction, namely, let
$k_{1i,s}(\vec{x},u)=\rho_s(u)(\frac12\cdot\vec{\varepsilon}_2+\frac{1}{4}\cdot
k_{1i,s-1}(\vec{x}))$ for $i<s\leq p$. In the end we obtain:
$$k_{1i}=k_{1i,p}: D^2_1\times D^2_i\times (S^1_2\times\cdots\hat S_i^1
\times\cdots \times S^1_{p})\hookrightarrow \RR^{p+2}.$$ From the
construction, we see that $k_{1i}$ can be extended a bit as an
embedding:
$$k_{1i}:B^4\times
(S^1_2\times\cdots\times\hat{S}^1_i\times\cdots\times S^1_{p})\hookrightarrow
\RR^{p+2}.$$
\end{proof}

\begin{cor}\label{ki} For each $i=1,\cdots,p$, $\imath_p$ also extends
as an embedding:
$$k_i:D^3\times(S^1_1\times\cdots\times \hat{S^1_i}\times\cdots\times
S^1_{p})\hookrightarrow \RR^{p+2},$$ where $S^1_1\times\cdots\times
S^1_p\subset D^3\times(S^1_1\times\cdots\times
\hat{S^1_i}\times\cdots\times S^1_{p})$ is an unknotted embedding
for $S^1_i$ in ${\rm Int}(D^3)$, and the identity on other
factors.\end{cor}

\begin{proof} Clearly the inclusion $S^1_1\times S^1_i\subset B^4\subset\CC\times\CC$ in Lemma \ref{k1i} extends as an
embedding $S^1\times D^3\cong(S^1_1\times D^1)\times B^2\subset
B^4\subset\CC\times\CC$ where $S^1_1\times D^1$ is a tubular
neighborhood of $S^1_1\subset\CC$ and $B^2$ is the $2$-dimensional
disk centered at the origin with radius $\frac{3}{2}$. Now $k_i$
may be defined as $k_{1i}$ composed with the latter embedding.
\end{proof}

We introduce the notion of unknotted embeddings and their types.

\begin{definition}\label{unknotted} A smooth embedding $\imath:T^p\hookrightarrow \RR^{p+2}$ is
called \emph{unknotted} if there is a compactly-supported self-diffeomorphism
$g:\RR^{p+2}\to\RR^{p+2}$ of $\RR^{p+2}$ such that $\imath$ and $g\circ\imath_p$
have the same image, i.e. $\imath(T^p)=g\circ\imath_p(T^p)$.
\end{definition}

\begin{remark}\label{rmkknot} Since the orientation-preserving
diffeomorphism group $\Diff_+(\RR^n)$
deformation-retracts to ${\rm SO}(n)$, and hence that $\pi_0\Diff_+(\RR^n)$ is trivial (\cite{St}),
clearly our definition of unknottedness
agrees with the more common notion that
$\imath(T^p)$ and $\imath_p(T^p)$ are equivalently knotted is
if there is a diffeotopy of $\RR^{p+2}$ taking $\imath(T^p)$ to $\imath_p(T^p)$.
\end{remark}

\begin{definition}\label{types}
Two unknotted embeddings $\imath_0,\imath_1:T^p\hookrightarrow \RR^{p+2}$ are
called of the same \emph{type} if they are the same up to a
self-diffeomorphism of $\RR^{p+2}$, namely there is a diffeomorphism
$h:\RR^{p+2}\to \RR^{p+2}$ such that $h\circ\imath_0=\imath_1$.
This is an equivalence relation, and the equivalent classes are
called \emph{types}. The type of $\imath_p$ is called the
\emph{standard} type.
For any $\tau\in\Aut(T^p)$, $\tau$ defines a \emph{modular
transformation} on types, namely $[\imath]\mapsto[\imath\circ\tau]$.
A \emph{modular} type is obtained by a modular transformation of the
standard type.
\end{definition}

\begin{lemma}\label{imageAndType} For any unknotted embedding and any type, there is an unknotted embedding
with the same image and of that type.\end{lemma}
\begin{proof} Let $\imath_0$ be the embedding, and $[\imath_1]$ be the type. By Definition \ref{types},
there is some $\RR^{p+2}$-self-diffeomorphism $h_1$ such that
$h_1\circ\imath_0(T^p)=\imath_1(T^p)$. Let $\tau=\imath_1^{-1}\circ
h_1\circ\imath_0:T^p\to T^p$, then
$h_1^{-1}\circ\imath_1=\imath_0\circ\tau^{-1}$. Thus
$\imath_0\circ\tau^{-1}$ has the same image as $\imath_0$, and the
same type as $[\imath_1]$.\end{proof}

Related is the notion of extendable automorphisms.

\begin{definition}\label{esg} Let $\imath:T^p\hookrightarrow\RR^{p+2}$ be a smooth embedding.
An automorphism $\tau\in\Aut(T^p)$ is said to be \emph{extendable} over $\imath$ if there is
a compactly-supported self-diffeomorphism of $\RR^{p+2}$ which commutes with $\tau$ via $\imath$. The subgroup of $\Aut(T^p)$
consisting of extendable automorphisms will be denoted as $E_\imath$.
\end{definition}

Note $E_{\imath\circ\tau}=\tau^{-1}E_{\imath}\tau$ for any smooth embedding $\imath$ and any $\tau\in\Aut(T^p)$. It is also clear that modular-type embeddings are in natural one-to-one correspondence with the
right cosets of $E_{\imath_p}$ in $\Aut(T^p)$, so Theorem \ref{main-extend} may be rephrased
as:

\begin{cor}\label{modTypeNumber} For $p\geq 1$, there are at most $2^p-1$ modular types of unknotted embeddings $\imath:T^p\hookrightarrow\RR^{p+2}$.\end{cor}

In the rest of this section, we prove Theorem \ref{main-extend}.

Fix a product structure
$T^p=S^1_1\times\cdots\times S^1_p$ as in Example \ref{stdmodel}, then
 $\Aut(T^p)$ is identified with $\SL(p,\ZZ)$ ($p\geq 2$). Denote:
$$R_{ij}=I+E_{ij},\quad Q_{ij}=R_{ij}^{-1}R_{ji}R_{ij}^{-1},$$ where $i\neq j$, and $I$ is the
identity matrix and $E_{ij}$ has $1$ for the $(i,j)$-entry and all
other entries $0$. Note that $R_{ij}$ is the full Dehn twist on
the sub-torus $S^1_i\times S^1_j$ along $S^1_i$, and $Q_{ij}$
trades the two factors of $S^1_i\times S^1_j$.

When $p=2$, there are two basic
extendable automorphisms for the embedding $T^2=S^1_1\times
S^1_2\subset \CC\times\CC\cong \RR^4$.

\begin{lemma}\label{R2Q} For the standard embedding $T^2=S^1_1\times S^1_2\subset B^4\subset\CC\times\CC\cong
\RR^4$, the following automorphisms can be extended as self-diffeomorphisms
of $B^4$ supported in the interior, (i.e. fixing an open neighborhood of $\partial B^4$):

(1) the twice full Dehn twist along each factor circle;

(2) trading two factors with their orientations preserved.
\end{lemma}

\begin{proof} (1) It suffices to prove for the first factor.
Consider $S_1^1\times S_2^1\subset (S_1^1\times D^1)\times D^2=
S_1^1\times D^3\subset \RR^4$, where $S_1^1\times D^1$ is a tubular
neighborhood of $S_1^1$ in the first  $\CC$, and $D^2$ is the disk
bounded by $S_2^1$ in the second $\CC$, $S_1^1\times D^3$ is a
tubular neighborhood of $S_1^1$ in ${\RR}^4$, $\partial (S^1_1
\times D^3)=S^1_1\times S^2$, and $*\times  S_2^1$ is the equator of
$*\times S^2$.

The Dehn $2$-twist $\tau :S_1^1\times S_2^1\to S_1^1\times S_2^1$ is
$(x,y)\mapsto (x,x^2y)$. The map $x\mapsto x^2$, considered as a map
from $S^1$ to $\SO(2)$, is of degree $2$. Thus the map $x\mapsto
x^2$, considered as a map $g:S^1\to \SO(3)$, is homotopic to a
constant map since $\pi_1(\SO(3))\cong \ZZ_2$. We may extend the map
$\tau :S_1^1\times S_2^1\to S_1^1\times S_2^1$ to a map $\tilde\tau
: S_1^1\times S^2\to
 S_1^1\times S^2$, defined by $\tau (x,y)=(x,g(x)y)$. Since
$g:S^1\to \SO(3)$ is homotopic to a constant map, $\tilde \tau$ is
diffeotopic to the identity. By the isotopy extension theorem (cf. \cite[Theorem 1.3]{Hi}),
$\tilde \tau$ can be extended to a self-diffeomorphism of $B^4$
supported in the interior.

(2) First extend $T^2=S^1_1\times S^1_2\to S^1_1\times S^1_2$ as
$f:\RR^4=\CC\times\CC\to\CC\times\CC$, by
$(z,w)\mapsto(\bar{w},z)$. $f$ is an orientation-reversing
diffeomorphism. To adjust to get an orientation-preserving one,
pick a self-diffeomorphism $h:\RR^4\to\RR^4$,
supported in the interior of $B^4\subset\RR^4$, such that $h(T^2)$ lies on the
subspace $\RR^3_{2,3,4}$. Let $r_1:\RR^4\to\RR^4$ be the
reflection with respect to $\RR^3_{2,3,4}$, then $r_1$ is
orientation-reversing, and is the identity restricted to $h(T^2)$.
Now $f_1=(h^{-1}r_1h)\circ f$ is orientation preserving, extending
the described automorphism on $T^2$. Furthermore, since
$f_1(z,w)=(-w,z)$ when $|w|^2+|z|^2>4-\varepsilon$, where
$\varepsilon>0$ is sufficiently small, we may adjust $f_1$ near
the boundary of $B^4$ to get an $f_2$ which is supported in the interior of $B^4$.
\end{proof}

From this observation we have the following lemma for $p\geq2$.

\begin{lemma}\label{pR2Q} $R^2_{1i},Q_{1i}\in E_{\imath_p}$, for $1<i\leq p$.\end{lemma}

\begin{proof} Let $\tau$ be either $R^2_{1i}$ or $Q_{1i}$.
From Lemma \ref{R2Q}, $\tau|_{S^1_1\times S^1_i}:S^1_1\times
S^1_i\to S^1_1\times S^1_i$ extends as $\bar{\tau}:B^4\to B^4$
which is the identity near the boundary. Therefore by Lemma
\ref{k1i},
$$\bar{\tau}\times\id:B^4\times
(S^1_2\times\cdots\times\hat{S}^1_i\times\cdots\times S^1_{p})\to
B^4\times (S^1_2\times\cdots\times\hat{S}^1_i\times\cdots\times
S^1_{p})$$ induces a self-diffeomorphism of $k_{1i}(B^4\times
(S^1_2\times\cdots\times\hat{S}^1_i\times\cdots\times S^1_{p}))$
which is the identity near the boundary. The latter further
extends to a diffeomorphism of $\RR^{p+2}$ by the identity outside
the image of $k_{1i}$.
\end{proof}

These automorphisms are not enough for generating
$E_{\imath_p}$ when $p\geq3$. Extra extendable ones come from the
following geometric construction.

\begin{lemma}\label{RR}
$R_{1p}R_{ip}\in E_{\imath_p}$, for $1<i<p$.
\end{lemma}

\begin{proof} According to Lemma \ref{k1i}, we will first extend
$\eta=R_{1p}R_{ip}$ on $T^p$ over $B^4\times
(S^1_2\times\cdots\times\hat{S}^1_i\times\cdots\times S^1_{p})$,
identity near the boundary, then extend it further as a
self-diffeomorphism of $\RR^{p+2}$ via $k_{1i}$ by the identity
outside the image. Since $\eta=\tau\times \id$ as from
$(S^1_1\times S^1_i\times S^1_p)\times
(S^1_2\times\cdots\times\hat{S}^1_i\times\cdots\times S^1_{p-1})$
to itself, essentially one must extend:
$$\tau:S^1_1\times S^1_i\times S^1_p\to S^1_1\times S^1_i\times
S^1_p,$$ as $B^4\times S^1_p\to B^4\times S^1_p$. It is easy to see
the matrix of $\tau$ is $
\left(\begin{array}{ccc}1&0&1\\0&1&1\\0&0&1
\end{array}\right),$
the $(1,i,p)$ minor of $R_{1p}R_{ip}$. It follows that each column
sum of $R_{1p}R_{ip}$ is odd. Rewrite $u_p$ as
$w$, then we have $\tau((u_1,u_i),w)=(\mu_{w}(u_1,u_i),w)$, where
$\mu_{w}(u_1,u_i)=(wu_1,wu_i)$ is an action of $S^1_p$ on $T^2$ by
flowing along the diagonal-slope direction.

We extend $\tau$ as follows. First, via the inclusion
$T^2=S^1_1\times S^1_i\subset S^3$ (in fact, the sphere centered
at the origin with radius $\sqrt{2}$), the diagonal-slope fibration
on $T^2$ extends as the Hopf fibration on $S^3$, and $\mu_{w}$
also extends as $\tilde{\mu}_{w}:S^3\to S^3$ by flowing along the
fiber loops. Thus $\tau$ extends as $\tilde{\tau}:S^3\times
S^1_p\to S^3\times S^1_p$,
$\tilde{\tau}(x,w)=(\tilde{\mu}_{w}(x),w)$.

On the other hand, $\tilde{\mu}_w$ can be regarded as
$S^1_p\to{\rm Diff}_+(S^3)$ which is given by a Lie group left
multiplication, regarding $S^1_p$ as a subgroup of $S^3$. Thus
$\tilde{\mu}_w$ extends as $S^3\to {\rm Diff}_+(S^3)$ by the Lie
group left multiplication. This implies that $\tilde{\mu}_w$ is
homotopic to the constant identity in $\pi_1{\rm Diff}_+(S^3)$
(since any circle is null-homotopic in the 3-sphere). Thus
$\tilde{\tau}$ is diffeotopic to the identity. By the isotopy
extension theorem, $\tilde{\tau}$ can be extended to a
diffeomorphism $\bar{\tau}$ of $B^4\times S_p^1$, which is the
identity near the boundary.

Now $\bar{\eta}=\bar{\tau}\times\id$ is a self-diffeomorphism of
$(B^4\times S^1_p)\times
(S^1_2\times\cdots\times\hat{S}^1_i\times\cdots\times S^1_{p-1})$,
or diffeomorphically, of $B^4\times
(S^1_2\times\cdots\times\hat{S}^1_i\times\cdots\times S^1_{p})$,
being the identity near the boundary. Finally extend $\bar{\eta}$
to a diffeomorphism $\RR^{p+2}\to \RR^{p+2}$ via $k_{1i}$ with the
identity outside the image.
\end{proof}

We need an elementary matrix lemma below. The technical result
(2) is needed in Section \ref{Sec-realization}.

\begin{lemma}\label{mx} (1) The subgroup $G$ of $\SL(p,\ZZ)$
generated by $R^2_{1i}$, $Q_{1i}$ ($1<i\leq p$) and $R_{1p}R_{ip}$
($1<i<p$) consists of all
the matrices $U\in\SL(p,\ZZ)$ whose entry sum of each column is odd.

(2) For any $1\leq i\leq p$, any $U\in\SL(p,\ZZ)$ can be written
as $KJ$ such that $K$ is a word in $R^2_{1j}$, $Q_{1j}$, $1< j\leq
p$, and $J$ has the minor $J^*_{ii}=1$.
\end{lemma}

\begin{proof} (1) Because
 $R^2_{ij}=Q_{1i}R^2_{1j}Q^{-1}_{1i}$, $Q_{ij}=Q_{1i}Q_{1j}Q^{-1}_{1i}$,
 for $1< i\neq j\leq p$, we have $R^2_{ij},Q_{ij}\in G$ ($i\neq j$). Also
 $R_{1k}R_{jk}=Q_{kp}^{-1}R_{1p}R_{jp}Q_{kp}$ ($k\neq1,j,p$),
 and $R_{ik}R_{jk}=Q_{1i}R_{1k}R_{jk}Q^{-1}_{1i}$ ($i\neq 1,j,k$),
 we have $R_{ik}R_{jk}\in G$ ($i,j,k$ mutually different). Multiplying by $R^2_{ij}$ from the left of a matrix
 adds twice of the $j$-th row to the $i$-th, and by $R_{ik}R_{jk}$
 adds the $k$-th row to both the $i$-th and the $j$-th, and by
 $Q_{ij}$ switches the $i$-th row and $j$-th row up to a sign.
We claim that any $U\in\SL(p,\ZZ)$
 with odd column sums becomes diagonal in $\pm 1$'s under
finitely many such operations, by a Euclid type algorithm described
below.

Suppose $U=(u_{ij})_{1\leq i,j\leq p}$. In the first column there
are an odd number of odd entries, say $u_{i_1,1},\cdots,
u_{i_{2m+1},1}$. Adding the $i_1$-th row simultaneously to the
$i_2$-th, $\cdots$, $i_{2m+1}$-th rows (i.e. multiplying $U$ from
the left by $R_{i_2i_1}R_{i_3i_1}\cdots
R_{i_{2m}i_1}R_{i_{2m+1}i_1}$) if necessary, we may assume only
one entry is odd, and switching that row with the first row if
necessary, we may further assume $u_{11}$ is the only odd entry in
the first column. Now there is some nonzero entry with a minimum
absolute value, say $u_{r1}$. As $U\in\SL(p,\ZZ)$, $u_{r1}$ cannot
divide all other entries in this column unless $u_{r1}=\pm1$, so
if $u_{r1}\neq \pm 1$, there must be some $u_{s1}$ with
$|u_{s1}|>|u_{r1}|$ and $u_{r1}$ cannot divide $u_{s1}$. By adding
(or subtracting) an even times of the $r$-th row to the $s$-th,
the $u_{s1}$ becomes $u_{s1}'$ such that
$-|u_{r1}|<u_{s1}'<|u_{r1}|$. Because this operation does not
change parity, we may repeat this process until we get a matrix
with $u_{11}=\pm1$, and other entries in the first column being
even. Then add (or subtract) several even times of the first row
to other rows, we obtain a matrix $U''$ with the entries in the
first column being zero except $u_{11}=\pm1$. Apply the process
recursively on the $(p-1)\times (p-1)$-submatrix
$U''_{11}=(u''_{ij})_{2\leq i,j\leq p}$, and use $u_{22}=\pm1$ to
kill other even entries in the second column, and so on. In the
end $U$ becomes a diagonal matrix with $\pm1$'s on the diagonal.

Finally, the claim means that any $U$ with odd column sums can be
written as $U=KD$ where $K$ is a word in $R^2_{ij}$, $Q_{ij}$,
$R_{ik}R_{jk}$, and $D$ is diagonal in $\pm 1$'s. Moreover, $D$
must have even $-1$'s on the diagonal since the determinant is
$1$, for example, at the $i_1,\cdots, i_{2m}$ place, $0\leq 2m\leq
p$, then $D=Q^2_{i_1i_2}\cdots Q^2_{i_{2m-1}i_{2m}}$. Therefore
$U=KD\in G$.

To see any element $U\in G$ has odd column sums, note that this
condition is the same as saying that $X\bar{U}=X$ where
$X$ is the row vector $(1,\cdots,1)\in\ZZ_2^{p}$, where $\bar{U}$ is the modulo $2$
reduction of $U$. As all the $R^2_{1i}$, $Q_{1i}$ ($1<i\leq p$)
and $R_{1p}R_{ip}$ ($1<i<p$) fix $X$, so does $G$. Therefore
$G$ consists of elements in $\SL(p,\ZZ)$ with odd column sums.

(2) Instead of doing row operations, do the first step of the
algorithm by column operations on the $i$-th row of $U^{-1}$ to
make the $(i,r)$-th entry $\pm 1$. Switch the $i$-th and $r$-th
column, and multiply $Q^2_{ir}$ if necessary to make the
$(i,i)$-th entry $+1$. In other words, the $(i,i)$-th entry of
$U^{-1}K$ is $1$ for some word $K$ in $R^2_{1i}$, $Q_{1i}$. Then
let $J=(U^{-1}K)^{-1}$.
\end{proof}

\begin{proof}[Proof of Theorem \ref{main-extend}.] By
Lemma \ref{pR2Q} and Lemma \ref{RR}, $G\leq E_{\imath_p}$, where $G$ is the subgroup
of $\SL(p,\ZZ)$ as in Lemma \ref{mx}. By Lemma
\ref{mx} (1), $G$ is the stabilizer of the row vector $(1,\cdots,1)\in\ZZ_2^p$
under the right action of $\SL(p,\ZZ)$ on the row vector space $\ZZ_2^p$.
Note that with the chosen product structure of $T^p$,
this action is naturally identified with the action of $\Aut(T^p)$ on $H^1(T^p;\ZZ_2)$,
thus $E_{\imath_p}$ contains the stabilizer of a nontrivial element in $H^1(T^p;\ZZ_2)$.
It follows that the
$[\Aut(T^p):E_{\imath_p}]\leq 2^p-1$ since $\Aut(T^p)$ acts invariantly and transitively on
the subset of nontrivial elements of $H^1(T^p;\ZZ_2)$.
\end{proof}

\section{Realizing DE attractors of $T^p$ in codimension $2$}\label{Sec-realization}

We prove Theorem \ref{main-realize-T^p} in this section. To fix the
notation, we write $(u_1,\cdots,u_p)$ for the coordinate of a point $u$ in
$T^p=S^1_1\times\cdots\times S^1_p$,
and write $(u,z)$ for the coordinate of a point in
$T^p\times D^2$. We use $T^{p-1}_i$ to denote
$S^1_1\times\cdots\times \hat{S^1_i}\times\cdots\times S^1_{p}$.

With the fixed product structure, expanding maps of $T^p$ may be identified
with expanding endomorphisms, i.e. which are represented by
$p\times p$ integral matrices whose eigenvalues all have absolute values strictly greater
than $1$. We identify expanding maps and automorphisms of $T^p$ with
 their matrices, and often write them as $\phi_A$,
$\tau_U$, etc.

Given an expanding map $\phi:T^p\to T^p$, we wish to use an
unknotted embedding to realize an attractor derived from $\phi$. It
is not hard to lift it as a hyperbolic bundle embedding first.
Note that the normal bundles of orientable codimension $2$
submanifold of $\RR^{p+2}$ must have trivial Euler classes, and hence are trivial bundles
(cf. \cite[Corollary 11.4]{MS}), we should consider self-embeddings $e$
of trivial disk bundles which lifts an expanding map $\phi$.
For simplicity, we also assume from now on that
$\phi$ has positive degree, possibly after passing to $\phi^2$.

\begin{prop} \label{bundle embedding} Let $\phi:T^p\to T^p$ be any orientation-preserving
expanding map.
There is a hyperbolic bundle embedding $e:T^p\times D^2\hookrightarrow
T^p\times D^2$ lifted from $\phi$.\end{prop}

\begin{proof} Let $A$ be the matrix of $\phi$. Recall that a Smith normal form ([Ne] Theorem II.9)
 of a positive-determinant integral matrix $A$ is a decomposition $A=U \Delta V$, where $U,V\in\SL(p,\ZZ)$,
and $\Delta$ is a diagonal matrix ${\rm diag}(\delta_1,\cdots,\delta_p)$,
with positive integral diagonal entries $\delta_i$, such that $\delta_i$
divides $\delta_{i+1}$ for each $1\leq i\leq p-1$. Let $\Delta_i= {\rm diag} (1, \cdots,
\delta_i, \cdots, 1)$, then: $$A=U \Delta_1 \cdots \Delta_p V.$$ We
first lift each factor to a bundle embedding $e_U, e_V,
e_{\Delta_i}: T^p\times D^2\hookrightarrow T^p\times D^2$.

To lift $\tau_U$, define, for example,
$e_U(u,z)=(\tau_U(u),u_1^{m_1}\cdots u_p^{m_p} z)$ for chosen
integers $m_1,\cdots,m_p$. Similarly lift $\tau_V$ to $e_V$. To
lift $\phi_{\Delta_i}$, first pick a hyperbolic bundle embedding
${b}_{\delta_i}:S^1_i\times D^2\hookrightarrow S^1_i\times D^2$ such that
${b}_{\delta_i}(u_i,z)=(u_i^{\delta_i},\bar b_i(u_i,z))$ sends the
solid torus into itself as a connected thickened closed braid with
winding number $\delta_i$, shrinking evenly on the disk direction.
Then define:
$$e_{\Delta_i}(u,z)=(\phi_{\Delta_i}(u),\bar b_i(u_i,z))=({b}_{\delta_i}\times \text{id}_{T^{p-1}_{i}})(u,z).$$

Finally, take the composition $e=e_U\circ e_{\Delta_1}\circ \cdots
\circ e_{\Delta_p}\circ e_V$, and we obtain a hyperbolic bundle
embedding lifted from $\phi$.
\end{proof}

Although the lifting is far from unique, for our purpose of use we must pick a
topologically simple one, composing which does not make the embedding `knottier' or `more twisted'
in $\RR^{p+2}$.

\begin{figure}[htb]
\centering
\psfrag{a}[]{\scriptsize{$b_{\delta_i}(S^1_i\times \{1\})$}}
\psfrag{b}[]{\scriptsize{{$S^1_i\times \{1\}$}}}
\psfrag{c}[]{\scriptsize{glue up without twist}}
\psfrag{d}[]{\scriptsize{the $(\delta_i,1)$-cable} $b_{\delta_i}$}
\includegraphics[scale=0.7]{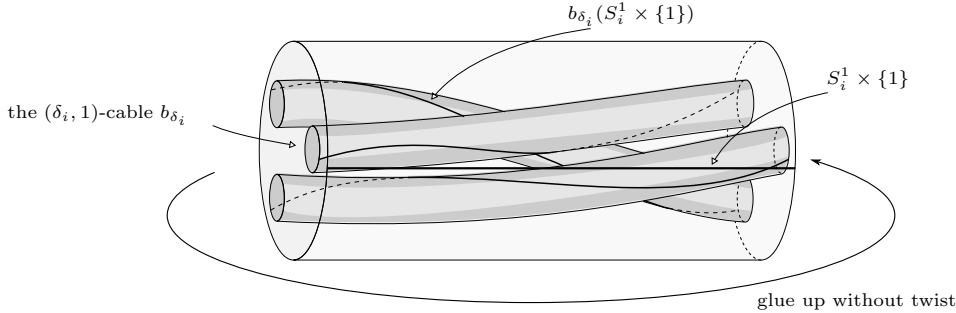}
\caption{The favorite lifting, illustrated as $\delta_i=3$.}\label{figFavLifting}
\end{figure}

\begin{example}[The favorite lifting]\label{favorite}  Define $e_U(u,z)=(\tau_U(u),z)$,
$e_V(u,z)=(\tau_V(u),z)$, and
$e_{\Delta_i}(u,z)={b}_{\delta_i}\times \text{id}_{T^{p-1}_{i}}$
with:
$${b}_{\delta_i}(u_i,z)=(u_i^{\delta_i}, \frac{1}{2}u_i+\frac{1}{\delta_i^2}u_i^{1-\delta_i}z).$$
The chosen braid $b _{\delta_i}$ is a $(\delta_i,1)$-cable with respect to 
the given trivialization of $S^1_{\delta_i}\times D^2$. It is presented in Figure \ref{figFavLifting} 
as $\delta_i=3$,
where the framing change is indicated by $b_{\delta_i}(S^1_i\times\{1\}$ and its image.

This specific choice of $b_{\delta_i}$ gives us a favorite lifting
$e=e_U\circ e_{\Delta_1}\circ \cdots \circ e_{\Delta_p}\circ e_V$
of $\phi$.
\end{example}

An embedding $\jmath:T^p\times D^2\hookrightarrow\RR^{p+2}$ is, in general,
understood by knowing its core restriction
$\imath=\jmath|_{T^p\times\{0\}}$, and the framing.
We need the following definition.

\begin{definition}\label{framing} $\jmath:T^p\times D^2\hookrightarrow \RR^{p+2}$ is called \emph{unknotted}, if its
core restriction $\imath$ is unknotted. The \emph{type} of $\jmath$
is the type of $\imath$. The embedding $\jmath$ is said to have
\emph{untwisted framing}, if $\jmath(T^p\times\{1\})$ is
null-homologous in the complement of the core
$\imath(T^p)$ in $\RR^{p+2}$.\end{definition}

\begin{prop}\label{utukJ} Unknotted embeddings with untwisted framing
are unique of their types, up to compactly-supported self-diffeomorphisms of $\RR^{p+2}$ .
Namely, for unknotted embeddings $\jmath_0,\jmath_1$ with untwisted
framing of the same type, there is a compactly-supported
self-diffeomorphism $h:\RR^{p+2}\to \RR^{p+2}$ of $\RR^{p+2}$ such that
$h\circ\jmath_0=\jmath_1$.\end{prop}

\begin{proof} By definition we may first find a compactly supported self-diffeomorphism $h$
of $\RR^{p+2}$ 
so that $\imath_1=h\circ\imath_0$, and in the smooth category one
may assume $h|:\jmath_0(T^p\times D^2)\to\jmath_1(T^p\times D^2)$ is conjugate to
the normal-bundle map, namely,
$$\jmath_1(u,z)=h(\jmath_0(u,g_u(z))),$$ for
$g_u\in\SO(2)$. Thus there is a continuous map $g:T^p\to\SO(2)$.
Note the set  $[T^p,\SO(2)]$ of homotopy classes of maps from $T^p$ to $\SO(2)$ is
in bijection to $H^1(T^p;\ZZ)\cong {\rm Hom}(H_1(T^p),\ZZ)$, as $\SO(2)\cong S^1$
is an Eilenberg-MacLane space $K(\ZZ,1)$. We may assume $g_u(z)=u_1^{m_1}\cdots u_p^{m_p}z$,
for some integers $m_1,\cdots,m_p$. In fact, $m_1,\cdots,m_p$ are
determined so that $g^*\xi=m_1\,[S^1_1]^*+\cdots+ m_p\,[S^1_p]^*$ in $H^1(T^p;\ZZ)$.
Here $\xi$ is the generator of $H^1(\SO(2))\cong\ZZ$ giving the natural orientation
of the normal bundle, and $([S^1_1]^*,\cdots,[S^1_p]^*)$
is the basis of $H^1(T^p)$ dual to the basis $([S^1_1],\cdots,[S^1_p])$
of $H_1(T^p)$.
Because $h(\jmath_0(u,1))$ is
null-homologous in $\RR^{p+2}-h(\imath_0(T^p))$ by the
untwisted-framing assumption of $\jmath_0$, it is not hard to check that the $p$-torus
defined by $h(\jmath_0(u,g_u(1))$, $u\in T^p$, represents
$g^*\xi$ under the Alexander duality
$H_p(\RR^{p+2}-\imath_1(T^p))\cong H^1(T^p)$. Because
$\jmath_1(u,1)=h(\jmath_0(u,g_u(1))$, this $p$-torus is also 
null-homologous in $\RR^{p+2}-\imath_1(T^p)$ by the
untwisted-framing assumption on $\jmath_1$. We see $g^*\xi=0$, so $m_i=0$
for $1\leq i\leq p$. Thus $g_u(z)=z$ for any $z\in D^2$, and $\jmath_1= h\circ \jmath_0$.
\end{proof}

Proposition \ref{utukJ} says a topologically `simple' embedding $\jmath$ is
determined by its type. We must show our favorite lifting $e$ is
topologically simple, namely composing $e$ preserves the untwisted-framing
property, the unknotted-ness. Moreover, we must show that any modular-type
embedding remains modular after composing
$e$, in order to use the finiteness result of Corollary \ref{modTypeNumber}. 
When $p=1$, there is no type issue, and the rest are the following well-known
facts in classical knot theory:

\begin{lemma}\label{simple}
(1) $b_{\delta_i}(S^1_i\times \{1\})$ is homological to $\delta_i$
times $S_i^1\times \{1\}$  in $S_i^1\times D^2 -
b_{\delta_i}(S^1_i\times \mathrm{Int}(D^2))$. Hence
if $\jmath: S^1_i\times D^2\hookrightarrow \RR^3$ has untwisted framing, then
$\jmath\circ b_{\delta_i}$ has untwisted framing.

(2) Furthermore if $\jmath$ in (1) is also unknotted, then
$b_{\delta_i}(S^1_i\times \{0\})$ is the $(\delta_i, 1)$-torus
knot, which is unknotted.
\end{lemma}

We prove the general case in Lemmas
\ref{stillut}, \ref{stilluk}.

\begin{lemma}\label{stillut} If $\jmath:T^p\times D^2\hookrightarrow\RR^{p+2}$ has untwisted
framing, then $\jmath\circ e$ also has untwisted framing.\end{lemma}

\begin{proof}

Obviously the composition with $e_U$, $e_V$
 preserves untwisted framing as they are identity on the fiber $D^2$.
  We claim $e_{\Delta_1},\cdots,e_{\Delta_p}$ preserves untwisted framing, then provided that
$\jmath$ has untwisted framing, so does $\jmath\circ e_U$, and hence
so does $(\jmath\circ e_U)\circ e_{\Delta_1}$, and so on. Finally
$\jmath\circ e$ also has untwisted framing.

Since $e_{\Delta_i}={b}_{\delta_i}\times \text{id}_{T^{p-1}_{i}}:
S^1_i\times D^2\times T^{p-1}_i\hookrightarrow S^1_i\times D^2\times
T^{p-1}_i$, by Lemma \ref{simple} (1), we have
$e_{\Delta_i}(T^p\times \{1\})$ is homological to $\delta_i$ times
$T^p\times \{1\}$ in $T^p\times D^2 -
e_{\Delta_i}(T^p\times\text{Int}(D^2))$. Thus if $\jmath:
T^p\times D^2\hookrightarrow \RR^{p+2}$ has untwisted framing, then $\jmath
\circ e_{\Delta_i}$ has untwisted framing.
\end{proof}

Showing that $e$ also preserves the unknotted-ness needs more effort. 
We first
investigate an essential case when $e$ is just $e_{\Delta_i}$, and
$\jmath$ is a special candidate of its type.

We will call $\jmath_p: T^p\times D^2 \hookrightarrow \RR^{p+2}$ {\it standard},
if $\jmath_p$ has untwisted framing, and $\jmath_p|_{T^p\times
\{0\}}=\imath_p$.

To help visualize how $\jmath\circ e(T^p\times D^2)$ unknots
itself, remember that a standard type $\jmath_p$ can
be written as $k_i\circ g_i$, for $i=1,...,p$,  where:
$$g_i:S^1_1\times\cdots\times S^1_p\times D^2\subset
(S^1_1\times\cdots\times \hat{S^1_i}\times\cdots\times
S^1_{p})\times D^3$$ is the identity on factors $S^1_j$, $j\neq
i$, and $S_i^1\times D^2\subset D^3$ is the thicken-up of the
circle lying on the equatorial disc of $D^3$, centered at the
origin with radius one half, (cf. Corollary \ref{ki}).

\begin{lemma}\label{stilluk-J} Suppose $1\leq i\leq p$, and
 $\jmath_J=\jmath_p\circ e_J$, where $J\in\SL(p,\ZZ)$
 satisfies the minor $J^*_{ii}=1$. Then
$\jmath_J\circ e_{\Delta_i}$ is also unknotted. Moreover, the type
of $\jmath_J\circ e_{\Delta_i}$ is modular.\end{lemma}

\begin{proof} Without loss of generality,
we may assume $i=1$. Denote $g_J=g_1\circ e_J$. The commutative
diagram below has included all the maps involved so far:
$$
 \xymatrix{
& & {T^{p-1}_1\times D^3}  \ar[dr]^{k_1} & \\
T^p\times D^2  \ar[d]^{\pi}\ar[r]^{e_{\Delta_1}} & T^p\times D^2
\ar[ur]^{g_J} \ar[d]^{\pi}
\ar[r]^{e_J} &  T^p\times D^2  \ar[d]^{\pi} \ar[u]^{g_1}  \ar[r]^{\jmath_p}  & {\RR^{p+2}}    \\
T^p  \ar[r]^{\phi_{\Delta_1}} & T^p   \ar[r]^{\tau_J} & T^p \ar[ur]^{\imath_p} }
$$

There is a basic picture for $p=2$ to
keep in mind. In this case, we are trying to unknot the core of $\jmath_J\circ e_{\Delta_1}(T^2\times D^2)$ 
in $\RR^4$. Let $K'|_1$ be a $(\delta_1,1)$-torus knot in $D^3$, whose
carrier torus is placed parallel to the $xy$-plane centered at the origin. $K'|_1$ is of course
an unknot. Let $r$ be an integer. Imagine an $S^1$-family of unknots $K'|_w$ in ${D}^3$, such that
for any $w\in S^1$, $K'|_w$ is obtained by rotating $K'|_1$ about the $z$-axis by
an angle $r\,{\rm arg}(w)$. We may simultaneously cap off all these knots by
picking a disk bounded by $K'|_1$ and rotate it about the $z$-axis 
so that at the time $w$ 
it is bounded by $K'|_w$. See Figure \ref{figCapRotate}. 
\begin{figure}[htb]
\centering
\psfrag{0}[]{\scriptsize{$D^3_0$}} \psfrag{1}[]{\scriptsize{$D^3_{\frac{\pi}{2}}$}}
\psfrag{2}[]{\scriptsize{$D^3_{\pi}$}}\psfrag{3}[]{\scriptsize{$D^3_{\frac{3\pi}2}$}}
\psfrag{w}[]{$w$}
\includegraphics[scale=0.7]{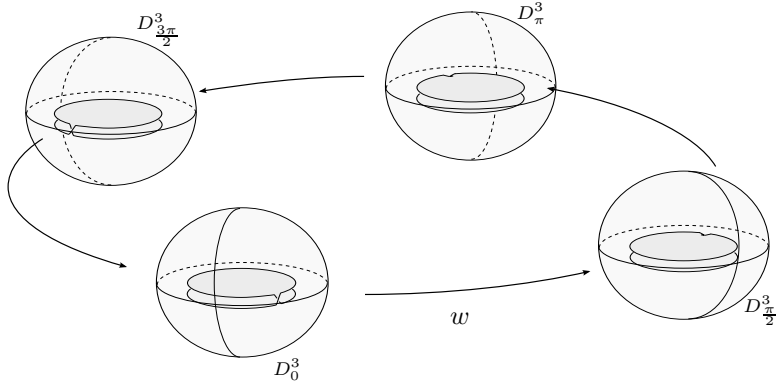}
\caption{An $S^1$-family of rotating $(\delta_1, 1)$-torus knots, illustrated as $\delta_1=2$.}\label{figCapRotate}
\end{figure}
This implies
that $K'=\bigcup_{w\in S^1}K'|_w$ is an unknotted torus in $S^1\times
D^3$ (in the sense that it is diffeotopic to $S^1\times S^1\subset
S^1\times D^3$). Therefore, if $S^1\times D^3$ is further
embedded in $\RR^4$, the image of $K'$ is also
an unknotted torus in $\RR^4$. As we will see below, $K'$ is exactly
$g_{J}\circ e_{\Delta_1}(T^2\times\{ 0\})$. 
When $p>2$, there is a similar picture, where we will have a $T^{p-1}_1$-family instead of just an $S^1$-family.
Another difference is that when $p=2$, the smooth 
mapping class group $\pi_0{\rm
Diff}_+(T^2)$ is isomorphic to $\SL(2,\ZZ)$ so the new embedding is automatically
modular, but in the general case, this is no longer true so we need to analyze more carefully.

Specifically, we wish to diffeotope the core
$K'=g_J\circ e_{\Delta_1}(T^p\times\{0\})$ back to
$K=g_1(T^p\times\{0\})$ within $T^{p-1}_1\times D^3$, before
including the latter into $\RR^{p+2}$ via $k_1$. The $(\delta_1,1)$-torus knot
$K'\subset T^{p-1}_1\times D^3$ may be viewed as a
$T_1^{p-1}$-family of (hopefully) unknotted loops in $D^3$.
One may regard $S^1_2,\cdots,S^1_p$ as independent clocks, and at
every `moment' $(u_2,\cdots,u_p)$ we see a loop in $D^3$. It turns
out that this is a $(\delta_1,1)$-torus knot rotating around a
fixed axis. Then we may simultaneously diffeotope the loops back to the
standard place in $D^3$. The key to reading out this picture is understanding the intersection
loci of $K'$ on fibers $D^3$, namely
$K'|_{(u_2,\cdots,u_p)}=K'\cap(\{(u_2,\cdots,u_p)\}\times D^3)$.

Denote $v=\phi(u)=\tau_J\circ\phi_{\Delta_1}(u)$, and
$\vec{\omega}(u_1)=u_{1x}\cdot\vec{\varepsilon}_1+u_{1y}\cdot\vec{\varepsilon}_2$
a rotating vector in $\RR^2_{1,2}$. Then:
$$g_1(u,z)=((u_2,\cdots,u_p),\ \frac{1}{2}\cdot\vec{\omega}(u_1)+\frac{1}{3}[z_x\cdot\vec{\omega}(u_1)+z_y\cdot\vec{\varepsilon}_3]).$$

\begin{figure}[htb]
\centering
\psfrag{b}[]{\scriptsize{$\bar{\varepsilon}_1$}}
\psfrag{a}[]{\scriptsize{$\bar{\varepsilon}_2$}}
\psfrag{c}[]{\scriptsize{$\bar{\varepsilon}_3$}}
\psfrag{e}[]{\scriptsize{$\frac{1}{3}(z_x\vec{\omega}(u_1)+z_y\bar{\varepsilon}_3)$}}
\psfrag{d}[]{\scriptsize{$\frac{1}{2}\vec{\omega}(u_1)$}}
\includegraphics[scale=.7]{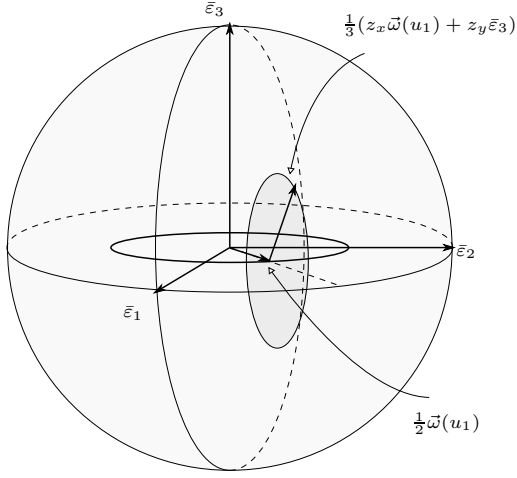}
\caption{The image of $(u_1,z)$ under $g_1:S_1^1\times D^2 \hookrightarrow D^3$.}\label{figMapBall}
\end{figure}

Composing $e_J(u,z)=(\tau_J(u),z)$ with
$e_{\Delta_1}(u,z)=(\phi_{\Delta_1}(u),\frac{1}{2}
u_1+\frac{1}{\delta_1^2}u_1^{1-\delta_1}z)$, we have:
$$g_J\circ
e_{\Delta_1}(u,0)=((v_2,\cdots,v_p),\ b(u)),$$ where $v_j$ means
the $S^1_j$-component of $v=\phi(u)$, and:
$$b(u)=\frac{1}{2}\cdot\vec{\omega}(v_1)+\frac{1}{6}[u_{1x}\cdot\vec{\omega}(v_1)+u_{1y}\cdot\vec{\varepsilon}_3].$$

Because $J^*_{11}=1$ is the $(1,1)$-th entry of $J^{-1}$,
$u_1^{\delta_1}=v_1v_2^{-r_2}\cdots v_p^{-r_p}$ for some integers
$r_2,\cdots,r_p$, so $v_1=u_1^{\delta_1}v_2^{r_2}\cdots
v_p^{r_p}$. Therefore at the `moment' $(v_2,\cdots,v_p)\in
T^{p-1}_1$, $K'|_{(v_2,\cdots,v_p)}$ is a $(\delta_1,1)$-torus
knot in $D^3$ defined by $b(u)$, and as $v_j=e^{i\theta_j}$, $K'|_{(v_2,\cdots,v_p)}$ rotates about the
$\vec{\varepsilon}_3$-axis by an angle $r_j\theta_j$.

Note that a $(\delta_1,1)$-torus knot in $D^3$ is unknotted, so
there is a diffeotopy $h_t:D^3\to D^3$ supported in the interior, such that $h_0=\id_{D^3}$ and
$h_1(K'|_{(1,\cdots,1)})=K|_{(1,\cdots,1)}$. To unknot $K'$
simultaneously on fibers, define $\rho:T^{p-1}_1\to \SO(3)$, with
$\rho {(v_2,\cdots,v_p)}$ being the rotation about the
$\vec{\varepsilon}_3$-axis by an angle
$r_2\theta_2+\cdots+r_p\theta_p$, where $v_j=e^{i\theta_j}$. The
`unknotting' diffeotopy may be defined as:
$$H_t:T^{p-1}_1\times D^3\to T^{p-1}_1\times D^3$$
with:
$$H_t((v_2,\cdots,v_p),\vec{x})=((v_2,\cdots,v_p),\rho (v_2,\cdots,v_p)\circ
h_t\circ (\rho (v_2,\cdots,v_p))^{-1}(\vec{x})).$$

Since $H_t$ is supported in the interior, when we embed
$T^{p-1}_1\times D^3$ into $\RR^{p+2}$ by $k_1$, $H_t$ induces a
diffeotopy on $k_1(T^{p-1}_1\times
D^3)$ supported in the interior, which extends as a diffeotopy of $\RR^{p+2}$. This
diffeotopy takes $\jmath_J\circ e_{\Delta_1}(T^p\times \{ 0\})$
back to $\jmath_p(T^p\times \{ 0\})$, so the former is unknotted
too.

To see the `moreover' part, note that both $g_1,H_1\circ g_J\circ
e_{\Delta_1}$ embed $T^p\times\{ 0\}$ into $T_1^{p-1}\times D^3$
with the same image. Let $B$ denote the matrix obtained from $J$
by multiplying each entry of the first column by $\delta_1$ and
then replacing the first row by $(1,0,\cdots,0)$. Since
$J^*_{11}=1$, $B\in\SL(p,\ZZ)$. By comparing: $$H_1\circ g_J\circ
e_{\Delta_1}(u,0)=(v_2,\cdots,v_p,\rho (v_2,\cdots,v_p)\circ
h_1\circ (\rho (v_2,\cdots,v_p))^{-1}\circ b(u)),$$ with:
$$g_1\circ e_B(u,0)=(v_2,\cdots,v_p,\frac{1}{2}\vec{\omega}(u_1)),$$
we see that the self-diffeomorphism $(H_1\circ g_J\circ
e_{\Delta_1})\circ (g_1\circ e_B)^{-1}: T^p\to T^p$ can be written
as $F:T^{p-1}_1\times S^1\to T^{p-1}_1\times S^1$,
such that $F((v_2,\cdots,v_p),z)=((v_2,\cdots,v_p),f(v_2,\cdots,v_p)(z))$,
where $f:T^{p-1}_1\to {\rm Diff}_+(S^1)$. Because ${\rm
Diff}_+(S^1)\simeq \SO(2)$, $f$ is homotopic to
$f_0:T^{p-1}_1\to \SO(2)$, where
$f_0(v_2,\cdots,v_p)=v_2^{m_2}\cdots v_p^{m_p}$ for some integers
$m_2,\cdots,m_p$, and $F$ is diffeotopic to the
automorphism $F_0:T_1^{p-1}\times S^1\to T_1^{p-1}\times S^1$,
$F_0((v_2,\cdots,v_p),z)=((v_2,\cdots,v_p),v_2^{m_2}\cdots
v_p^{m_p}z)$. Let $M$ denote the matrix obtained from the $p\times
p$ identity matrix by replacing the first row by
$(1,m_2,\cdots,m_p)$. Then the type of $\jmath_J\circ
e_{\Delta_1}$ is a modular transformation of the standard type by
$MB$.
\end{proof}

The general case that composing $e$ preserves the unknotted-ness and remains modular is
as follows.

\begin{lemma}\label{stilluk} If $\jmath$ is unknotted of modular
type with untwisted framing, then $\jmath\circ e:T^p\times D^2\hookrightarrow
\RR^{p+2}$ is also unknotted of modular type.\end{lemma}

\begin{proof} Clearly
$\hat\jmath=\jmath\circ e_U$ is unknotted with modular type and
untwisted framing. Note $\hat\jmath\circ e_{\Delta_1}$ is
unknotted if and only if so is $h\circ\hat\jmath\circ
e_{\Delta_1}$ for any $\RR^{p+2}$ self-diffeomorphism $h$. By
Proposition \ref{utukJ}, the unknotted-ness of $\hat\jmath\circ
e_{\Delta_1}$ depends only on the type of $\hat\jmath$. Since
$\hat\jmath$ is of modular type, suppose $\jmath_J=\jmath_p\circ
e_J$ has the same type as $\hat\jmath$. Moreover, $J$ may be
picked so that $J^*_{11}=1$ by Lemma \ref{mx}(2). Then by Lemma
\ref{stilluk-J}, $\jmath_J\circ e_{\Delta_1}$ remains unknotted
with modular type. Therefore $\hat\jmath\circ e_{\Delta_1}$ is
unknotted with modular type. By Lemma \ref{stillut},
$\hat\jmath\circ e_{\Delta_1}$ has untwisted framing.

Repeat this argument so $(\hat\jmath\circ e_{\Delta_1})\circ
e_{\Delta_2}$ is unknotted with modular type and untwisted
framing, and so on we see that $\hat\jmath\circ e_{\Delta_1}\cdots
e_{\Delta_p}$ is also unknotted with modular type. Finally
$(\jmath \circ e_U\circ e_{\Delta_1}\circ \cdots \circ
e_{\Delta_p})\circ e_V$ is also unknotted since it has the same
image of the core as $\jmath\circ e_U\circ e_{\Delta_1}\circ
\cdots \circ e_{\Delta_p}$, and the type change is modular. We
conclude that $\jmath\circ e$ is still unknotted of modular type.
\end{proof}

\begin{proof}[Proof of Theorem \ref{main-realize-T^p}]
Pick an unknotted embedding $\jmath:T^p\times D^2\hookrightarrow \RR^{p+2}$ of
modular type with untwisted framing. (We can in fact require the
core image to be any unknotted $T^p$ in $\RR^{p+2}$, by
Lemma \ref{imageAndType}.) By Lemma \ref{stillut} and Lemma
\ref{stilluk}, $\jmath\circ e^i$ ($i\geq 0$) are all unknotted of
modular type with untwisted framing. By Corollary \ref{modTypeNumber}, at least two of
$\jmath\circ e^i$ ($0\leq i\leq 2^p-1$) are of the same type.
Suppose $\jmath\circ e^k$ and $\jmath\circ e^l$ ($0\leq k<l\leq
2^p-1$) are of the same type. Pick the embedding
$\hat\jmath=\jmath\circ e^k:T^p\times D^2\hookrightarrow \RR^{p+2}$ instead of
$\jmath$, and let $d=l-k$. There is a compactly-supported
self-diffeomorphism $h$ of $\RR^{p+2}$ such that
$h\circ\hat\jmath=\hat\jmath\circ e^d$ by Proposition \ref{utukJ}.
This is to say, $e^d$ can be realized by an embedding
$\hat\jmath:T^p\times D^2\hookrightarrow\RR^{p+2}$ with extension $h$.

Therefore, $$\Lambda=\bigcap_{i=0}^\infty e^i(T^p\times
D^2)=\bigcap_{i=0}^\infty e^{di}(T^p\times D^2),$$ embeds into
$\RR^{p+2}$ by $\hat\jmath$ as an attractor of $h$, so we have
realized an expanding attractor derived from $\phi$.
\end{proof}

\bibliographystyle{amsplain}

\begin{thebibliography}{}
\setlength{\itemsep}{0ex} \addcontentsline{toc}{section}{Reference}

\bibitem[Bo]{Bo} H. Bothe,  {\it The ambient structure of expanding attractors. II.
Solenoids in $3$-manifolds.} Math. Nachr. \textbf{112} (1983), 69--102.

\bibitem[DLWY]{DLWY} F. Ding, Y. Liu,  S. Wang,  J. Yao, {\it Spin structure and codimension-two
homeomorphism extension.} Preprint, 2010. \texttt{arXiv:0910.4949}.

\bibitem[DPWY]{DPWY} F. Ding, J. Pan, S. Wang, J. Yao, {\it Manifolds with $\Omega(f)$ a union of DE
attractors are rational homology spheres}. Ergodic Theory Dynam. Systems \textbf{30} (2010), no. 5, 1399--1417.

\bibitem[DSS]{DSS} K. Dekimpe, M. Sadowski, A Szczepa\'nski, {\it Spin structures on flat manifolds}. (English
summary). Monatsh. Math. \textbf{148} (2006), no. 4, 283--296.

\bibitem [ES]{ES} D. Epstein, M. Shub, {\it Expanding endomorphisms on flat manifolds}. Topol. \textbf{7}
(1968), 139--141.

\bibitem[Gr]{Gr} M. Gromov, \textit{Groups of polynomial growth and expanding
maps.} Inst. Hautes \'Etudes Sci. Publ. Math. \textbf{53} (1981), 53--73.

\bibitem[Ha]{Ha} A. Haefliger, \textit{Plongements diff\'erentiables dans le domaine stable}.
(French). Comment. Math. Helv. \textbf{37} (1962/1963), 155--176.

\bibitem[Hi]{Hi} M. Hirsch, \textit{Differential Topology}. Graduate Text in Mathematics, Vol. 33.
Springer-Verlag, New York-Heidelberg, 1976.


\bibitem[JNW]{JNW} B. Jiang, Y.  Ni,  S. Wang, {\it 3-manifolds that admit knotted
solenoids as attractors.} Trans. Amer. Math. Soc. \textbf{356} (2004), no.
11, 4371--4382.

\bibitem[JWZ]{JWZ} B. Jiang, S. Wang, H. Zheng, {\it No embeddings of solenoids
into surfaces.} Proc. Amer. Math. Soc. \textbf{136} (2008), no. 10,
3697--3700.

\bibitem[MS]{MS} J. Milnor, J. Stasheff, {\it Characteristic classes.} Annals
of Mathematics Studies 76. Princeton University Press, Princeton, NJ;
University of Tokyo Press, Tokyo, 1974.

\bibitem[Ne]{Ne} M. Newman, \textit{Integral Matrices}.
Pure and Applied Mathematics, Vol. 45. Academic Press, New York-London, 1972.


\bibitem[Sm]{Sm} S. Smale, \textit{Differentiable dynamical systems},
    Bull. Amer. Math. Soc. \textbf{73} (1967), 747--817.

\bibitem[St]{St} T.~E. Stewart, \textit{On groups of diffeomorphisms.}
Proc. Amer. Math. Soc. \textbf{11} (1960), 559--563.

\bibitem[Wu]{Wu} W.~T. Wu, \textit{On the isotopy of
$C^{r}$-manifolds of dimension $n$ in euclidean $(2n+1)$-space}.
Sci. Record (N.S.) \textbf{2} (1958), 271--275.

\end{thebibliography}

\bigskip
\textsc{Peking University,
Beijing, 100871, P.~R. China.}

\textit{E-mail address:} \texttt{dingfan@math.pku.edu.cn}

\bigskip
\textsc{University of California,
Berkeley, CA 94720, USA.}

\textit{E-mail address:} \texttt{yliu@math.berkeley.edu}

\bigskip
\textsc{Peking University,
Beijing, 100871, P.~R. China.}

\textit{E-mail address:} \texttt{wangsc@math.pku.edu.cn}

\bigskip
\textsc{University of California,
Berkeley, CA 94720, USA.}

\textit{E-mail address:} \texttt{jiangangyao@gmail.com}

\end{document}